\documentclass[sn-mathphys, pdflatex]{sn-jnl}

\jyear{2022}%
\usepackage{algpseudocode}
\usepackage{caption}
\usepackage{subcaption}
\theoremstyle{thmstyleone}%
\newtheorem{theorem}{Theorem}
\theoremstyle{thmstyletwo}%

\theoremstyle{thmstylethree}%
\usepackage{amssymb}
\usepackage{array}

\begin{document}

\title[Bragin, M. A. and Tucker, E. L.: Surrogate ``Level-Based" LR for MILP]{Surrogate ``Level-Based" Lagrangian Relaxation for Mixed-Integer Linear Programming}


\author[1]{\fnm{Mikhail A.} \sur{Bragin$^*$}}\email{mikhail.bragin@uconn.edu
}

\author[2]{\fnm{Emily L.} \sur{Tucker}}\email{etucke3@clemson.edu}

\affil[1]{\orgdiv{Department of Electrical \& Computer Engineering}, \orgname{University of Connecticut}, \orgaddress{\street{371 Fairfield Way, U-4157}, \city{Storrs}, \postcode{06269}, \state{CT}, \country{USA}}}

 \affil[2]{\orgdiv{Department of Industrial Engineering}, \orgname{Clemson University}, \orgaddress{\street{271 Freeman Hall}, \city{Clemson}, \state{SC}, \postcode{29634}, \country{USA}}}



\abstract{Mixed-Integer Linear Programming (MILP) plays an important role across a range of scientific disciplines and within areas of strategic importance to society. The MILP problems, however, suffer from \textit{combinatorial complexity}. Because of integer decision variables, as the problem size increases, the number of possible solutions increases \textit{super-linearly} thereby leading to a drastic increase in the computational effort. To efficiently solve MILP problems, a ``price-based'' decomposition and coordination approach is developed to exploit 1. the super-linear reduction of complexity upon the decomposition and 2. the geometric convergence potential inherent to Polyak’s stepsizing formula for the fastest coordination possible to obtain near-optimal solutions in a computationally efficient manner. Unlike all previous methods to set stepsizes heuristically by adjusting hyperparameters, the key novel way to obtain stepsizes is purely decision-based: a novel ``auxiliary'' constraint satisfaction problem is solved, from which the appropriate stepsizes are inferred. Testing results for large-scale Generalized Assignment Problems (GAP) demonstrate that for the majority of instances, certifiably optimal solutions are obtained. For stochastic job-shop scheduling as well as for pharmaceutical scheduling, computational results demonstrate the two orders of magnitude speedup as compared to Branch-and-Cut (B\&C). The new method has a major impact on the efficient resolution of complex Mixed-Integer Programming (MIP) problems arising within a variety of scientific fields.}




\maketitle
\vspace{-6mm}
\section{Introduction}\label{sec1}

Mixed-Integer Linear Programming (MILP) plays an important role across a range of scientific disciplines such as mathematics, operations research, engineering, and computer science as well as within a range of areas of strategic importance to society such as biology \cite{Huang16,Yang21}, healthcare \cite{Kayvanfar2021, Hachicha16}, humanitarian applications \cite{smalley2015optimized, hamdan2020robust, Ahani21, kamyabniya2021robust}, manufacturing \cite{Liu21a, hong2019admission, balogh2022milp, oztop2022metaheuristics}, pharmacy \cite{Kopanos10, Stefannson11, Zhu2018, Ge21}, power and energy systems \cite{schill2017start, Yonghong20, Li20}, transportation and logistics \cite{archetti2021optimization, reddy2022multi} and many others. 

The associated systems are created by interconnecting $I$ smaller subsystems, each having its own objective and a set of constraints. The subsystem interconnection is modeled through the use of \textit{system-wide coupling} constraints. Accordingly, the MILP problems are frequently formulated in terms of cost components associated with each subsystem with the corresponding objective functions being additive as such:
\begin{flalign}
& \min_{(x,y) := \left\{x_i,y_i\right\}_{i=1}^I} \Bigg\{\sum_{i=1}^I \left((c_i^x)^T x_i + (c_i^y)^T y_i\right)  \Bigg\}.  \label{eq1}
\end{flalign}
Furthermore, 
coupling constraints 
are additive in terms of $I$ subsystems: 
\begin{flalign}
& s.t. \;\; \sum_{i=1}^I A_i^x x_i + \sum_{i=1}^I A_i^y y_i - b = 0, \;\; \left\{x_i,y_i\right\} \in \mathcal{F}_i, i = 1, \dots, I.  \label{eq2}
\end{flalign}
The \textit{primal} problem \eqref{eq1}-\eqref{eq2} is assumed to be feasible and the feasible region $\mathcal{F} = \prod_{i=1}^I  \mathcal{F}_i$ with $ \mathcal{F}_i \subset \mathbb{Z}^{n_i^x} \times \mathbb{R}^{n_i^y}$
is assumed to be bounded and finite. 
The MILP problems modeling the above systems are referred to as \textit{separable}. Because of the discrete decisions, however, MILP problems are known to be NP-hard and are prone to the curse of \textit{combinatorial complexity}. As the size of a problem increases, the associated number of combinations of possible solutions (hence the term ``combinatorial'') increases super-linearly (e.g., exponentially) thereby making problems of practical sizes difficult to solve to optimality; even near-optimal solutions are frequently difficult to obtain.

A beacon of hope to resolve combinatorial difficulties lies through the exploitation of separability through the \textit{dual} ``price-based'' decomposition and coordination Lagrangian Relaxation technique. After the relaxation of coupling constraints \eqref{eq2}, the coordination of subproblems amounts to the maximization of a concave non-smooth dual function:
\vspace{-2mm}
\begin{flalign}
& \max_{\lambda} \{q(\lambda): \lambda \in \mathbb{R}^m\},  \label{eq3}
\end{flalign}
where
\vspace{-2mm}
\begin{flalign}
& q(\lambda) \! = \!\min_{(x,y)} \Bigg\{L(x,y,\lambda), \label{eq4}   \{x_i,y_i\} \in \mathcal{F}_i, i = 1,\dots,I \Bigg\}. 
\end{flalign}
Here $L(x,y,\lambda) \equiv \sum_{i=1}^I (c_i^x)^T x_i \!+\! \sum_{i=1}^I (c_i^y)^T y_i \!+\! \lambda^T \cdot \left(\sum_{i=1}^I A_i^x x_i \!+\! \sum_{i=1}^I A_i^y y_i \! -\! b\right)$ is the Lagrangian function. The Lagrangian multipliers $\lambda$ (``dual'' variables) are the decision variables with respect to the dual problem \eqref{eq3}, and it is assumed that the set of optimal solutions is not empty. The minimization within \eqref{eq4} with respect to $\{x,y\}$ is referred to as the ``relaxed problem.'' 

While the sizes of the primal and the relaxed problems are the same in terms of the number of discrete variables, the main advantage of Lagrangian Relaxation is the exploitation of the reduction of the combinatorial complexity upon decomposition into subproblems. Accordingly, the number of discrete decision variables within the primal problem is $n = \sum_{i=1}^I n_i^x$, so the worst-case complexity of solving the primal problems is $O(e^{\sum_{i=1}^I n_i^x})$. By the same token, the worst-case complexity required to solve the following subproblem
\begin{flalign}
& \min_{x_i,y_i} \left\{(c_i^x)^T x_i + (c_i^y)^T y_i + \lambda^T \cdot \left(A_i^x x_i + A_i^y y_i \right),\left\{x_i,y_i\right\} \in \mathcal{F}_i\right\}, 
\label{eq8}
\end{flalign}
is $O(e^{n_i^x})$. The decomposition ``reverses'' the combinatorial complexity thereby exponentially reducing the effort. The decomposition, therefore, offers a viable potential to improve the operations of existing systems as well as to scale up the size of the systems to support their efficient operations.



\begin{figure}[!ht]
  \centering
    \includegraphics[trim=210 480 180 70,  width=0.4\linewidth, scale=0.5]{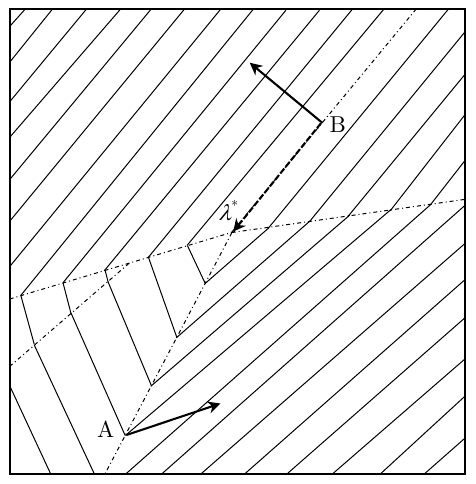}
    \caption{An example of a dual function demonstrating difficulties faced by subgradient methods. Solid lines denote the level curves, dash-dotted lines denote the ridges of the dual function whereby the usual gradients are not defined (possible subgradient directions at points A and B are denoted by solid arrows), and the direction from point B toward optimal multipliers is denoted by a dashed line.}
    \label{fig_ex1}
\end{figure}

While decomposition efficiently reduces the combinatorial complexity, the coordination aspect of the method to efficiently obtain the optimal ``prices'' (Lagrangian multipliers) has been the subject of an intense research debate for decades because of the fundamental difficulties of non-smooth optimization. Namely, because of the presence of integer variables $x$, the dual function \eqref{eq3} is non-smooth comprised of flat convex polygonal facets (each corresponding to a particular solution to the relaxed problem within \eqref{eq4}) intersecting at linear ridges along which the dual function $q(\lambda)$ is non-differentiable; in particular, $q(\lambda)$ is not differentiable at $\lambda^{*}$ thereby ruling out the possibility of using necessary and sufficient conditions for the extremum. As a result of the non-differentiability of $q(\lambda)$,
subgradient multiplier-updating directions, however, are non-ascending directions thereby leading to a decrease of dual values; subgradient directions may also change drastically thereby resulting in zigzagging of Lagrangian multipliers (see Figure \ref{fig_ex1} for illustrations) and slow convergence as a result.


Traditional methods to maximize $q(\lambda)$ rely upon iterative updates of Lagrangian multipliers by taking a series of steps $s^k$  along subgradient $g(x^k,y^k)$ directions as:
\vspace{-2mm}
\begin{flalign}
& \lambda^{k+1} = \lambda^k + s^k \cdot g(x^k,y^k),   \label{eq5}
\end{flalign}
where $\{x^k,y^k\} \equiv \{x^k_i,y^k_i\}_{i=1}^I$ is a an optimal solution to the relaxed problem \eqref{eq4} with multipliers equal to $\lambda^k.$ 
Within the Lagrangian Relaxation framework, subgradients are defined as levels of constraint violations $g(x^k,y^k) \equiv \sum_{i=1}^I A_i^x x_i^k + \sum_{i=1}^I A_i^y y_i^k- b.$ 
Inequality constraints $\sum_{i=1}^I A_i^x x_i + \sum_{i=1}^I A_i^y y_i \leq b$, if present, can be handled by converting into equality constraints by introducing non-negative real-valued slack variables $z$ such that $\sum_{i=1}^I A_i x_i  + \sum_{i=1}^I A_i^y y_i + z = b.$ The multipliers are subsequently projected onto the positive orthant delineated by restrictions $\lambda \geq 0.$

Because of the lack of differentiability of $q(\lambda)$, notably, at the optimum $\lambda^{*},$ the stepsize selection plays an important role to guarantee convergence to the optimum as well as for the success of the overall Lagrangian Relaxation methodology for solving MILP problems. 

One of the earlier papers on the optimization of non-smooth convex functions, with $q(\lambda)$ being its member, though irrespective of Lagrangian Relaxation, is Polyak's seminal work \cite{Polyak69}. 
Intending to achieve the geometric (also referred to as ``linear'') rate of convergence so that $\|\lambda^k - \lambda^*\|$ is monotonically decreasing, Polyak proposed the stepsizing formula, which in terms of the problem under consideration takes the following form:  
\vspace{-2mm}
\begin{flalign}
& 0 < s^k < \gamma \cdot \frac{q(\lambda^{*}) - q(\lambda^k)}{\left\|g(x^k,y^k)\right\|^2}, \gamma < 2.  \label{eq7}
\end{flalign}
Within \eqref{eq7} and thereafter in the paper
the standard Euclidean norm is used unless stated otherwise. 

Subgradient directions, however, 1. are generally difficult to obtain computationally when the number of subproblems \eqref{eq8} to be solved is large, and 2. change drastically thereby resulting in zigzagging of Lagrangian multipliers and slow convergence. Moreover, 3. stepsizes \eqref{eq7} cannot be set due to the lack of knowledge about the optimal dual value $q(\lambda^*)$. 

To overcome the first two of the difficulties above, the Surrogate Subgradient Method was developed by \cite{Zhao99} whereby the exact optimality of the relaxed problem (or even subproblems) is not required. As long as the following ``surrogate optimality condition'' is satisfied: 
\vspace{-2mm}
\begin{flalign}
& L(\tilde{x}^k,\tilde{y}^k,\lambda^k) < L(\tilde{x}^{k-1},\tilde{y}^{k-1},\lambda^k)  \label{eq9}
\end{flalign}
the multipliers can be updated by using the following formula
\vspace{-2mm}
\begin{flalign}
& 0 < s^k < \gamma \cdot \frac{q(\lambda^{*}) - L(\tilde{x}^k,\tilde{y}^k,\lambda^k)}{\left\|g(\tilde{x}^k,\tilde{y}^k)\right\|^2}, \;\; \gamma < 1,  \label{eq10}
\end{flalign}
and convergence to $\lambda^{*}$ is guaranteed. Here ``tilde'' is used to distinguish optimal solutions $\{x^k\}$ to the relaxed problem, from the solutions $\{\tilde{x}^k\}$ that satisfy the ``surrogate optimality condition'' \eqref{eq9}. Unlike that in Polyak's formula, parameter $\gamma$ is less than 1 to guarantee that $q(\lambda^{*}) > L(\tilde{x}^k,\tilde{y}^k,\lambda^k)$ so that the stepsizing formula \eqref{eq10} is well-defined, as proved in \cite[Proposition 3.1, p. 703]{Zhao99}.
Once $\{\tilde{x}^k,\tilde{y}^k\}$ are obtained, multipliers are updated by using the same formula as in \eqref{eq5} with stepsizes from \eqref{eq10} and ``surrogate subgradient'' multiplier-updating directions $g(\tilde{x}^k,\tilde{y}^k)$ used in place of subgradient directions $g(x^k,y^k)$. 
Besides reducing the computational effort owing to \eqref{eq9}, the concomitant reduction of multiplier zigzagging has also been observed. The main difficulty is the lack of knowledge about $q(\lambda^*)$. As a result, the geometric/linear convergence of the method (or any convergence at all) is highly questionable. Nevertheless, the underlying geometric convergence principle behind the formula \eqref{eq9} is promising and will be exploited in Section \ref{Results}. 

One of the first attempts to overcome the difficulty associated with the unavailability of the optimal [dual] value is the subgradient-level method developed by \cite{Goffin99} by adaptively adjusting a ``level'' estimate based on the detection of ``sufficient descent'' of the [dual] function and ``oscillation'' of [dual] solutions. In a nutshell, a ``level'' estimate is set as $q^{k}_{lev} = q^{k_j}_{rec} + \delta_j$ with $q^{k}_{rec}$ being the best dual value (``record objective value'') obtained up to an iteration $k,$ and $\delta_j$ is an adjustable parameter with $j$ denoting the $j^{th}$ update of $q^{k}_{lev}.$ Every time oscillations of multipliers are detected, $\delta_j$ is reduced by half. In doing so, stepsizes appropriately decrease,  $q^{k}_{lev}$ increases (for maximization problems) and the process continues until $\delta_j \rightarrow 0$ and  $q^{k}_{lev} \rightarrow q(\lambda^*).$ 

To improve convergence, rather than updating all the multipliers ``at once,'' within the incremental subgradient methods \cite{Nedic01}, multipliers are updated ``incrementally.'' 
Convergence results of the subgradient-level method (\cite{Goffin99}) have been extended for the incremental subgradient methods and proved. 

Within the Surrogate Lagrangian Relaxation (SLR) Method \cite{SLR}, the computational effort is reduced along the lines of the Surrogate Subgradient Method \cite{Zhao99} discussed above, that is, by solving one of a few subproblems at a time. To guarantee convergence, within SLR, distances between multipliers at consecutive iterations are required to decrease through a specially-constructed contraction mapping until convergence. Within \cite[Figs. 3-5,7; pp. 195-199]{SLR}, the SLR method showed the advantage against the above-mentioned subgradient-level method \cite{Goffin99} and the incremental subgradient methods \cite{Nedic01b, Nedic01} for non-smooth optimization. Unlike the methods of \cite{Nedic01b, Nedic01}, the SLR method does not require obtaining dual values to set stepsizes, which further reduces the effort. 
Aiming to simultaneously guarantee convergence while ensuring fast reduction of constraint violations and preserving the linearity of the original MILP problem, the Surrogate Absolute-Value Lagrangian Relaxation (SAVLR) method \cite{Bragin19} was developed to penalize constraint violations by using $l_1$ ``absolute-value'' penalty terms. The above methods are reviewed in more detail in Section \ref{Supp}.

Because of the presence of the integer variables, there is the so-called the \textit{duality gap}, which means that even at convergence, $q(\lambda^*)$ is generally less than the optimal cost of the original problem \eqref{eq1}-\eqref{eq2}. To obtain a feasible solution to \eqref{eq1}-\eqref{eq2}, the subproblem solutions when put together may not satisfy all the relaxed constraints. Therefore, to solve corresponding MILP problems, heuristics are inevitable and are used to perturb subproblem solutions. The important remark here is that the closer the multipliers are to the optimum, generally, the closer the subproblem solutions are to the global optimum of the original problem, and the easier it is to obtain feasible solutions through heuristics. Therefore, having fast convergence in the dual space to maximize the dual function \eqref{eq3} is of paramount importance for the overall success of the method. Specific heuristics will be discussed at the end of Section \ref{Results}.
\section{Results}\label{Results}
\subsection{Surrogate ``Level-Based'' Lagrangian Relaxation} \label{SLBLR}

In this subsection, a novel Surrogate ``Level-Based'' Lagrangian Relaxation (SLBLR) method is developed to determine “level” estimates of $q(\lambda^*)$ within the Polyak's stepsizing formula for fast convergence of multipliers when optimizing the dual function \eqref{eq3}. Since the knowledge of $q(\lambda^*)$ is generally unavailable, over-estimates of the optimal dual value, if used in place of $q(\lambda^*)$ within the formula \eqref{eq10}, may lead to the oscillation of multipliers and to the divergence. Rather than using heuristic ``oscillation detection'' of multipliers used to adjust ``level'' values \cite{Goffin99}, the key of SLBLR is the decision-based ``divergence detection'' of multipliers based on a novel auxiliary ``multiplier-divergence-detection” constraint satisfaction problem. 

\noindent \textbf{``Multiplier-Divergence-Detection'' Problem to Obtain the Estimate of $q(\lambda^*).$} The premise behind the multiplier-divergence detection is the rendition of the result due \cite[Theorem 4.1, p. 706]{Zhao99}:

\begin{theorem}
Under the stepsizing formula 
\begin{flalign} 
& s^k < \gamma \cdot \frac{q(\lambda^{*}) - L(\tilde{x}^k,\tilde{y}^k,\lambda^k)}{\|g(\tilde{x}^k,\tilde{y}^k)\|^2}, \gamma < 1, \label{eq16}
\end{flalign}
such that $\{\tilde{x}^k,\tilde{y}^k\}$ satisfy
\begin{flalign} 
& L(\tilde{x}^k,\tilde{y}^k,\lambda^k) \leq L(\tilde{x}^{k-1},\tilde{y}^{k-1},\lambda^k),\label{eq17}
\end{flalign}
the multipliers move closer to optimal multipliers $\lambda^{*}$ iteration by iteration:  
\begin{flalign} 
& \|\lambda^{*}-\lambda^{k+1}\| < \|\lambda^{*}-\lambda^{k}\|.\label{eq18}
\end{flalign}
\end{theorem}

The following Corollary and Theorem 2 are the main key results of this paper.

\noindent \textbf{Corollary 1}
\textit{If} 
\begin{flalign} 
& \|\lambda^{*}-\lambda^{k+1}\| \geq \|\lambda^{*}-\lambda^{k}\|, \label{eq19}
\end{flalign}
\textit{then} 
\begin{flalign} 
& s^k \geq \gamma \cdot \frac{q(\lambda^{*}) - L(\tilde{x}^k,\tilde{y}^k,\lambda^k)}{\|g(\tilde{x}^k,\tilde{y}^k)\|^2}.\label{eq20}
\end{flalign}

\begin{theorem}
If the following auxiliary ``multiplier-divergence-detection'' feasibility problem (with $\lambda$ being a continuous decision variable: $\lambda \in \mathbb{R}^m$)
\begin{flalign} 
    \begin{cases}
      \|\lambda-\lambda^{k_j+1}\| \leq \|\lambda-\lambda^{k_j}\|,\\
      \|\lambda-\lambda^{k_j+2}\| \leq 
\|\lambda-\lambda^{k_j+1}\|,\\
...\\
\|\lambda-\lambda^{k_j+n_j}\| \leq 
\|\lambda-\lambda^{k_j+n_j-1}\|, \label{eq23}
    \end{cases}
\end{flalign}
admits no feasible solution with respect to $\lambda$ for some $k_j$ and $n_j$, then $\exists \; \kappa \in [k_j,k_j+n_j]$ such that 
\begin{flalign} 
& s^\kappa \geq \gamma \cdot \frac{q(\lambda^{*}) - L(\tilde{x}^\kappa,\tilde{y}^\kappa,\lambda^\kappa)}{\|g(\tilde{x}^\kappa,\tilde{y}^\kappa)\|^2}. \label{eq24}
\end{flalign}
\end{theorem}

\proof{} 
Assume the contrary: $\forall \kappa \in [k_j,k_j+n_j]$ the following holds:  
\begin{flalign} 
& s^\kappa < \gamma \cdot \frac{q(\lambda^{*}) - L(\tilde{x}^\kappa,\tilde{y}^\kappa,\lambda^\kappa)}{\|g(\tilde{x}^\kappa,\tilde{y}^\kappa)\|^2}. \label{eq25}
\end{flalign}
By Theorem 1, multipliers approach $\lambda^{*},$ therefore, the ``multiplier-divergence-detection'' problem admits at least one feasible solution - $\lambda^*.$ Contradiction. 
\endproof

\noindent From \eqref{eq24} it follows that $\exists \; \overline{q}_{\kappa,j}$ such that $\overline{q}_{\kappa,j} > q(\lambda^{*})$ and the following holds: 
\begin{flalign} 
& s^\kappa = \gamma \cdot \frac{\overline{q}_{\kappa,j} - L(\tilde{x}^\kappa,\tilde{y}^\kappa,\lambda^{\kappa})}{\|g(\tilde{x}^\kappa,\tilde{y}^\kappa)\|^2}. \label{eq26}
\end{flalign}
The equation \eqref{eq26} can equivalently be rewritten as: 
\begin{flalign} 
& \overline{q}_{\kappa,j} = \frac{1}{\gamma} \cdot s^{\kappa} \cdot \|g(\tilde{x}^\kappa,\tilde{y}^\kappa)\|^2 + L(\tilde{x}^\kappa,\tilde{y}^\kappa,\lambda^\kappa). \label{eq27}
\end{flalign}
Therefore, 
\begin{flalign} 
& \overline{q}_{j} = \max_{\kappa \in [k_j,k_j+n_j]} \overline{q}_{\kappa,j} > q(\lambda^{*}). \label{eq28}
\end{flalign}

A brief yet important discussion is in order here. The overestimate $\overline{q}_{j}$ of the dual value $q(\lambda^*)$ is the sought-for ``level'' value after the $j^{th}$ update (the $j^{th}$ time the problem \eqref{eq23} is infeasible). Unlike previous methods, which require heuristic hyperparameter adjustments to set level values, within SLBLR, level values are obtained by using the decision-based principle per \eqref{eq23} precisely when divergence is detected without any guesswork. In a sense, SLBLR is hyperparameter-adjustment-free. Specifically, neither ``multiplier-divergence-detection'' problem \eqref{eq23}, nor the computations within \eqref{eq26}-\eqref{eq28} requires hyperparameter adjustment. Following \cite{Nedic01b}, the parameter $\gamma$ will be chosen as a fixed value $\gamma = \frac{1}{I}$, which is the inverse of the number of subproblems and will not require further adjustments. 

 Note that \eqref{eq23}  simplifies to an LP constraint satisfaction problem. For example, after squaring both sides of the first inequality $\|\lambda-\lambda^{k_j+1}\| \leq \|\lambda-\lambda^{k_j}\|$ within \eqref{eq23}, after using the binomial expansion, and canceling $\|\lambda-\lambda^{k_j}\|^2$ from both sides, the inequality simplifies to $2 \cdot (\lambda - \lambda^{k_j}) \cdot g(\tilde{x}^{k_j},\tilde{y}^{k_j}) \geq s^{k_j} \cdot \| g(\tilde{x}^{k_j},\tilde{y}^{k_j}) \|^2,$ which is linear in terms of $\lambda.$ 

To speed up convergence, a hyperparameter $\zeta < 1$ will be introduced to reduce stepsizes as follows: 
\begin{flalign} 
& s^{k} = \zeta \cdot \gamma \cdot \frac{\overline{q}_{j} - L(\tilde{x}^{k},\tilde{y}^{k},\lambda^{k})}{\|g(\tilde{x}^{k},\tilde{y}^{k})\|^2}, \zeta < 1. \label{eq29}
\end{flalign}
Subsequently, after iteration $k_{j+1}$ the problem \eqref{eq23} is sequentially solved again by adding one inequality per multiplier-updating iteration until iteration $k_{j+1}+n_{j+1}-1$ is reached for some $n_{j+1}$ so that \eqref{eq23} is infeasible. Then, stepsize is updated by using $\overline{q}_{j+1}$ per \eqref{eq29} and is used to update multipliers until the next time it is updated to $\overline{q}_{j+2}$ when the ``multiplier-divergence-detection'' problem is infeasible again, and the process repeats. Per \eqref{eq29}, SLBLR requires hyperparameter $\zeta$, yet, it is set before the algorithm is run and subsequently is not adjusted (see Numerical Testing Section \ref{Example:Ex_GAP} for empirical demonstration of the robustness of the method with respect to the choice of hyperparameter $\zeta$). 

To summarize the advantage of SLBLR, hyperparameter adjustment is not needed. The guesswork of when to adjust the level-value, and by how much is obviated – after \eqref{eq23} is infeasible, the level value is formulaically recalculated.

\noindent \textbf{On Improvement of Convergence.} To speed up the acceleration of the multiplier-divergence detection through the ``multiplier-divergence-detection'' problem, \eqref{eq23} is modified, albeit heuristically, in the following way: 
\begin{flalign} 
    \begin{cases}
      \|\lambda-\lambda^{k_j+1}\| \leq \sqrt{1-2 \cdot \nu \cdot s^{k_j}} \cdot \|\lambda-\lambda^{k_j}\|,\\
      \|\lambda-\lambda^{k_j+2}\| \leq 
\sqrt{1-2 \cdot \nu \cdot s^{k_j+1}} \cdot \|\lambda-\lambda^{k_j+1}\|,\\
...\\
\|\lambda-\lambda^{k_j+n_j}\| \leq 
\sqrt{1-2 \cdot \nu \cdot s^{k_j+n_j-1}} \cdot \|\lambda-\lambda^{k_j+n_j-1}\|. \label{eq31}
    \end{cases}
\end{flalign}
Unlike the problem \eqref{eq23}, the problem \eqref{eq31} no longer simplifies to an LP problem. Nevertheless, the system of inequalities delineates the convex region and can still be handled by commercial software. 

\noindent \textbf{Discussion of \eqref{eq31}.} Equation \eqref{eq31} is developed based on the following principles: 1. Rather than detecting divergence per \eqref{eq23}, convergence with a rate slower than $\sqrt{1-2 \cdot \nu \cdot s}$ is detected. This will lead to a faster adjustment of the level values. While the level value may no longer be guaranteed to be the upper bound to $q(\lambda^*)$, the merit of the above scheme will be empirically justified in the Numerical Testing Subsection \eqref{Testing}.
2. While the rate of convergence is unknown, in the “worst-case” scenario $\sqrt{1-2 \cdot \nu \cdot s}$ is upper bounded by 1 with $\nu = 0$, thereby reducing \eqref{eq31} to \eqref{eq23}. The estimation of $\sqrt{1-2 \cdot \nu \cdot s}$ is thus much easier than the previously used estimations of $q(\lambda^*)$ (as in Subgradient-Level-Based and Incremental approaches). 3. As the stepsize approaches zero the rate of convergence $\sqrt{1-2 \cdot \nu \cdot s}$ approaches one regardless of the value of $\nu$, once again reducing \eqref{eq31} to \eqref{eq23}. 
\vspace{2mm}

\noindent \textbf{Algorithm: Pseudocode.}  

\textbf{Input} $\lambda^0,$ $\gamma$, $\nu$, $\zeta$, $\overline{q}_0=+\infty$, $q^{max}=-\infty$
\begin{algorithmic}[1]
\While {$\frac{f(x^{feas},y^{feas}) - q(\lambda^k)}{f(x^{feas},y^{feas})} > \varepsilon_{gap}$} 
\State solve subproblem \eqref{eq8} s.t. \eqref{eq9}, 
\State calculate $g(\tilde{x}^{k},\tilde{y}^{k})$
\State calculate $L(\tilde{x}^{k},\tilde{y}^{k},\lambda^{k})$
\State calculate stepsizes per \eqref{eq29} as $s^{k} = \zeta \cdot \gamma \cdot \frac{\overline{q}_{j} - L(\tilde{x}^{k}\!,\tilde{y}^{k},\lambda^{k})}{\|g(\tilde{x}^{k}\!,\tilde{y}^{k})\|^2}$
\State update multipliers per \eqref{eq5} by using $g(\tilde{x}^k,\tilde{y}^k)$ as $\lambda^{k+1} = \lambda^k + s^k \cdot g(\tilde{x}^k,\tilde{y}^k)$
\If{$q^{max}\! < \!s^k \!\!\cdot\! \frac{\|g(\tilde{x}^{k}\!,\tilde{y}^{k})\|^2}{\gamma}\! +\! L(\tilde{x}^{k}\!,\tilde{y}^{k}\!,\lambda^{k})$} $q^{max}\! = \!s^k \!\cdot\! \frac{\|g(\tilde{x}^{k}\!,\tilde{y}^{k})\|^2}{\gamma} + L(\tilde{x}^{k}\!,\tilde{y}^{k}\!,\lambda^{k})$
\EndIf
\State $i \gets i+1$
\State $k \gets k+1$
\If{$i = I$} $i \gets 1$ 
\EndIf 
\If{\eqref{eq23} is infeasible} $\overline{q}_j = q^{max}, q^{max} = -\infty, j \gets j + 1$ 
\EndIf
\If{$\frac{\overline{q}_j - q(\lambda^k)}{\overline{q}_j} < \varepsilon$} search for feasible solutions $x^{feas},y^{feas}$ that satisfy 

$\! \! \! \! \!$ \eqref{eq2} to obtain a feasible cost 

$\! \! \! f(x^{feas},y^{feas}) \equiv \sum_{i=1}^I \left((c_i^x)^T x_i^{feas} + (c_i^y)^T y_i^{feas}\right)$
\EndIf
\EndWhile
\end{algorithmic}

There are three things to note here. 1. Steps in lines 15-16 are optional since other criteria can be used such as the number of iterations or the CPU time; 2. The value of $q(\lambda^k)$ is still needed (line 1) to obtain a valid lower bound. To obtain $q(\lambda^k)$, all subproblems are solved optimally for a given value of multipliers $\lambda^k$. The frequency of the search for the value $q(\lambda^k)$ is determined based on criteria as stated in point 1 above; 3. The search for feasible solutions is explained below.

\noindent \textbf{Search for Feasible Solutions.} Due to non-convexities caused by discrete variables, the relaxed constraints are generally not satisfied through coordination, even at convergence. Heuristics are thus inevitable, yet, they are the last step of the feasible-solution search procedure. Throughout all examples considered, following \cite{Bragin19} (as discussed in Section \ref{OnlineMethods}), $l_1$-absolute-value penalties penalizing constraint violations are considered. After the total constraint violation reaches a small threshold value, a few subproblem solutions obtained by the Lagrangian Relaxation method are perturbed, e.g., see heuristics within accompanying CPLEX codes within \cite{Bragin19} to automatically select which subproblem solutions are to be adjusted to eliminate the constraint violation to obtain a solution feasible with respect to the overall problem. 




\subsection{Numerical Testing}\label{Testing}

In this subsection, a series of examples are considered to illustrate different aspects of the SLBLR method.
%
In Example \ref{SmallEx}, a small example with known corresponding optimal Lagrangian multipliers is considered to test the new method as well as to compare how fast Lagrangian multipliers approach their optimal values as compared to Surrogate Lagrangian Relaxation \cite{SLR} and to Incremental Subgradient \cite{Nedic01} methods. 
In Example \ref{Example:Ex_GAP}, 
large-scale instances of generalized assignment problems (GAPs) of types D and E with 20, 40, and 80 machines and 1600 jobs from the OR-library (\href{https://www-or.amp.i.kyoto-u.ac.jp/members/yagiura/gap/}{https://www-or.amp.i.kyoto-u.ac.jp/members/yagiura/gap/}) are considered to demonstrate efficiency, scalability, robustness, and competitiveness of the method with respect to the best results available thus far in the literature. 
In Example \ref{JSS}, a stochastic version of a job-shop scheduling problem instance with 127 jobs and 19 machines based on \cite{Luh93} is tested. In Example \ref{pharma}, two instances of pharmaceutical scheduling with 30 and 60 product orders, 17 processing units, and 6 stages based on \cite{Kopanos10} are tested.  

For Examples \ref{SmallEx} and \ref{Example:Ex_GAP}, SLBLR is implemented within CPLEX 12.10 by using a Dell Precision laptop Intel(R) Xeon(R) E-2286M CPU @ 2.40GHz with 16 cores and installed memory (RAM) of 32.0 GB. For Examples \ref{JSS} and \ref{pharma}, SLBLR is implemented within CPLEX 12.10 by using a server Intel(R) Xeon(R) Gold 6248R CPU @ 3.00GHz with 48 cores and installed memory (RAM) of 192.0 GB.

\subsubsection{Demonstration of Convergence of Multipliers Based on a Small Example with Known Optimal Multipliers.} \label{SmallEx}
\begin{figure}[h]
  \centering
    \includegraphics[trim=10 0 0 0, width=1\linewidth, scale=0.3]{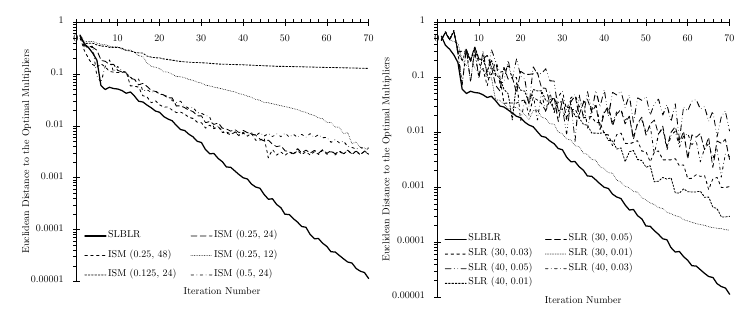}
    \caption{Results for Example \ref{SmallEx}: Comparison of SLBLR to 1. Incremental Subgradient Method (left) and 2. Surrogate Lagrangian Relaxation Method (right).}
    \label{Fig2}
\end{figure}

To demonstrate the convergence of multipliers, consider the following example (due \cite{Bragin20}): 
\begin{flalign}
& \min_{x_1, x_2, x_3, x_4, x_5, x_6} \left\{ x_1 + 2x_2 +3x_3 +x_4+2x_5+3x_6\right\}, \label{eq12c}
\end{flalign}
\begin{flalign}
& s.t. \; x_1 + 3x_2 +5x_3 +x_4+3x_5+5x_6 \geq 26, \\ & 2x_1 + 1.5x_2 +5x_3 +2x_4+0.5x_5+x_6 \geq 16. \label{eq13b}
\end{flalign}
As proved in \cite{Bragin20}, the optimal dual solutions are $\lambda_1^{*} = 0.6$ and $\lambda_2^{*} = 0.$ 
Inequality constraints are converted to equality constraints after introducing slack variables. 
In Figure \ref{Fig2}, the decrease of the corresponding distances from current multipliers to the optimal multipliers ($\|\lambda^k-\lambda^*\|$) is shown, and the SLBLR method is compared with the Incremental Subgradient method \cite{Nedic01} and the Surrogate Lagrangian Relaxation method \cite{SLR}. 

Within the SLBLR method, the equation \eqref{eq23} is used to detect divergence, and $\zeta = \frac{1}{2}$ is used to set stepsizes within \eqref{eq29}. In essence, only one hyperparameter was required, which has a quite simple explanation - ``when the stepsize is `too large,' cut the stepsize in half.'' As demonstrated in Figure \ref{Fig2}, the SLBLR method converges fast with $\|\lambda^k-\lambda^*\|$ decreasing roughly along a straight line on a log-scale graph suggesting that the rate of convergence is likely linear as expected. 

As for the Incremental Subgradient method, two hyperparameters are required: $R$ and $\delta$ (corresponding values used are shown in parentheses in the legend of Figure \ref{Fig2} (left)). A trial-and-error analysis indicated that ``acceptable'' values are $R=0.25$ and $\delta = 24.$ Increasing or decreasing $R$ to $0.5$ and $0.125$, respectively, do not lead to improvements. Likewise, increasing or decreasing $\delta$ to 48 and 12, respectively, do not lead to improvements as well. ``Plateau'' regions in the figure are caused by the fact that as stepsizes get smaller, a larger number of iterations is required for multipliers to travel the predetermined distance $R$; during these iterations, stepsizes are not updated and multipliers may oscillate around a neighborhood of the optimum without getting closer. While the above difficulty can be alleviated and convergence can be improved by hyperparameters $\tau$, $\beta$, and $R_l$ as reviewed in Supplementary Information Section \ref{Supp}, however, a larger number of hyperparameters would be required. 

As for the SLR method, several pairs of hyperparameters ($M$ and $r$) have been used as well (corresponding values used are shown in parentheses in the legend of Figure \ref{Fig2} (right)), yet, without exceeding the performance of the SLBLR method. 

Herein lies the advantage of the novel SLBLR method: the decision-based principle behind computing the ``level'' values. This is in contrast to the problem-dependent choice of hyperparameters $R$ and $\delta$ within the Subgradient-Level \cite{Goffin99} and Incremental Subgradient \cite{Nedic01} methods, and the choice of $M$ and $r$ within Surrogate Lagrangian Relaxation \cite{SLR, Bragin19} (see Sections \ref{sec1} and \ref{Supp} for more detail). 

Even after obtaining ``appropriate'' values of the aforementioned hyperparameters, the procedure that entails effort, results obtained by Surrogate Lagrangian Relaxation and the Incremental Subgradient Method do not match or beat those obtained by the SLBLR method. The specific reasons are 1. To obtain level values, heuristic adjustments of the ``level'' values are required (\cite{Goffin99}, \cite{Nedic01}) based on multiplier ``oscillation detection'' or ``significant descent'' (for minimization of non-smooth functions). However, these rules do not detect whether multipliers ``start diverging.'' Moreover, oscillation of multipliers is a natural phenomenon when optimizing non-smooth functions as discussed in Section \ref{sec1} since multipliers may zigzag/oscillate across ridges of the function, so the multiplier ``oscillation detection'' may not necessarily warrant the adjustment of level values. On the other hand,  multiplier ``oscillation'' is detected by checking whether multipliers traveled a (heuristically) predetermined distance $R$, hence, the divergence of multipliers can go undetected for a significant number of iterations (hence, the ``plateau'' regions shown in Figure \ref{Fig2} (left)), depending on the value of $R$. To the best of the authors' knowledge, the subgradient- and surrogate-subgradient-based methods using Polyak's (or Polyak-like) stepsizes with the intention of achieving the geometric/linear convergence rate either require $q(\lambda^*)$, which is unavailable, or require multipliers to travel infinite distance to guarantee convergence to the optimum $\lambda^*$ \cite{Goffin99}. 2. While SLR avoids the need to estimate $q(\lambda^*)$, the geometric/linear convergence is only possible outside of a neighborhood of $\lambda^*$ \cite[p. 187]{SLR}. Precisely for this reason, the convergence of multipliers within SLR with the corresponding stepsizing parameters $M=30$ and $r=0.01$ (as shown in Figure \ref{Fig2} (right)) appears to follow closely convergence within SLBLR up until iteration 50, after which the improvement tapers off.



\subsubsection{Generalized Assignment Problems.} \label{Example:Ex_GAP}

To demonstrate the computational capability of the new method as well as to determine appropriate values for key hyperparameters $\zeta$ and $\nu$ while using standard benchmark instances, large-scale instances of GAPs are considered (formulation is available in \ref{Formulation_GAP}). We consider 20, 40, and 80 machines with 1600 jobs (https://www-or.amp.i.kyoto-u.ac.jp/members/yagiura/gap/).

To determine values for $\zeta$ within \eqref{eq29} and $\nu$ within \eqref{eq31} to be used throughout the examples, several values are tested using GAP instance d201600. 

\begin{table}[!htbp]
\caption{Robustness Results for Instance d201600 with Respect to $\zeta$.} 
\centering 
\begin{tabular}{c c c c c} 
\hline\hline 
 $\zeta$ & Feasible  & Gap  & ``Auxiliary''  & Total   \\
 & Cost & (\%) & Time (sec) & Time (sec)  \\[0.5ex] 
\hline 

$1/1.25$  &     97827 &	 0.0059\%
 & 4.59 & 2904.02
 \\

$1/1.5$ &     \textbf{97825} & 0.0037\%
 & 17.10 & 1195.36
 \\

$1/2$ &     \textbf{97825} & 0.0048\%
 &	 88.59 &	2612.48 \\
 
$1/4$ &    97827 & 0.0059\%
 &	89.01 &	10235.50 \\

\hline 
\end{tabular}
\label{table0.5a} 
\end{table}

In Table \ref{table0.5a}, with fixed values of $\nu = 2$ and $s^0 = 0.02$, the best result (both in terms of the cost and the CPU time) is obtained with $\zeta = 1/1.5$. With the value of $\zeta = 1/4,$ the stepsize decreases ``too fast'' thereby leading to a larger number of iterations and a much-increased CPU time as a result. 

\begin{table} [!htbp]
\caption{Robustness Results for Instance d201600 with Respect to $\nu$.} 
\centering 
\begin{tabular}{c c c c c} 
\hline\hline 
 $\nu$ & Feasible  & Gap  & ``Auxiliary''  & Total   \\
 & Cost & (\%) & Time (sec) & Time (sec)  \\[0.5ex] 
\hline 

0.03125  &     97826 &	 0.0048\%
 & 93.79 & 2716.68 \\

0.125  &     \textbf{97825} & 0.0037\%
 & 33.62 & 1820.96
 \\

0.5 &     97826 & 0.0048\%
 &	 9.61 &	2444.46 \\
 
2 &     \textbf{97825} & 0.0037\%
 &	17.10 & 1195.36 \\

\hline 
\end{tabular}
\label{table0.5b} 
\end{table}

Likewise, in Table \ref{table0.5b} with fixed values of $\zeta = 1/1.5$ and $s^0 = 0.02$, it is demonstrated that the best result (both in terms of the cost and the CPU time) is obtained with $\nu = 2$. Empirical evidence here suggests that the method is stable for other values of $\nu.$ 

The robustness with respect to initial stepsizes ($s^0$) is tested and the results are demonstrated in Table \ref{table1} using the GAP type D instance with 20 machines and 1600 jobs. Multipliers are initialized by using LP dual solutions. 
The method's performance is appreciably stable for the given range of initial stepsizes used (Table \ref{table1}). 

\begin{table}[!htbp]
\caption{Robustness Results for Instance d201600 with Respect to Initial Stepsizes $s^0$.} 
\centering 
\begin{tabular}{c c c c c} 
\hline\hline 
 Initial & Feasible  & Gap  & ``Auxiliary''  & Total   \\
 Stepsize ($s^0$) & Cost & (\%) & Time (sec) & Time (sec)  \\[0.5ex] 
\hline 

0.0025  &     \textbf{97825} &	 0.0037\%
 & 123.71 & 2427.71 \\

0.005  &     \textbf{97825} & 0.0037\%
 &	6.84 & 1226.17 \\

0.01 &     97826 & 0.0048\%
 &	 6.96 & 2143.58 \\
 
0.02 &     \textbf{97825} & 0.0037\%
 &	17.10 & 1195.36 \\

0.04 &     97826 & 0.0048\%
 &	 19.21 & 1941.55
    \\ 
\hline 
\end{tabular}
\label{table1} 
\end{table}
SLBLR is robust with respect to initial multipliers $\lambda^0$ (Table \ref{table2}). For this purpose, the multipliers are initialized randomly by using the uniform distribution $U[90,110].$  For the testing, the initial stepsize $s^0=0.02$ was used. As evidenced from Table \ref{table2}, the method's performance is stable, exhibiting only a slight degradation of solution accuracy and an increase of the CPU time as compared to the case with multipliers initialized by using LP dual solutions. 

\begin{table}[!htbp]
\caption{Robustness Results for Instance d201600 with Respect to Initial Multipliers $\lambda^0$.} 
\centering 
\begin{tabular}{c c c c c c} 
\hline\hline 
 Case & Feasible  & Total Subproblem   & Feasible Solution & ``Auxiliary''  & Total   \\
 Number & Cost & Solving Time  & Search Time & Time & Time  \\
  &  &  (sec)  &  (sec) & (sec) & (sec)  \\[0.5ex] 
\hline 

1  &     \textbf{97825} &	1098.74	& 375.96 &	22.13 &	1496.84
    \\ 
2  &     97826 &	 1009.42 &	777.16 &	173.48 &	1960.07
    \\ 
3 &     97826 &	 2223.99 &	221.70 &	4.54 &	2450.24 
     \\ 
4 &     97826 &	 2333.55 &	402.41 &	4.08 &	2740.04 
     \\ 
5 &     97826 &	 1002.77 &	119.91 &	160.73 &	1283.42
    \\ 
\hline 
\end{tabular}
\label{table2} 
\end{table}






To test the robustness of the method across several large-scale GAP instances as well as scalability, six instances d201600, d401600, d801600, e201600, e401600, and e801600 are considered. 
SLBLR is compared with 
Depth-First Lagrangian Branch-and-Bound Method (DFLBnB) \cite{Posta12}, Column Generation \cite{Sadykov15}, and Very Large Scale Neighborhood Search (VLNS)
\cite{Haddadi19}, which to the best of the authors' knowledge are the best methods for at least one of the above instances. For completeness, a comparison against Surrogate Absolute-Value Lagrangian Relaxation (SAVLR) \cite{Bragin19}, which is an improved version of 
Surrogate Lagrangian Relaxation (SLR) \cite{SLR}, is also performed. The latter SLR method \cite{SLR} has been previously demonstrated to be advantageous against other non-smooth optimization methods as explained in Section \ref{OnlineMethods}.    
Table \ref{table4} presents feasible costs and times for each method. The advantage of SLBLR is the ability to obtain optimal results across a wider range of GAP instances as compared to other methods. Even though the comparison in terms of the CPU time is not entirely fair, feasible-cost-wise, SLBLR decisively beats previous methods. 
For the d201600 instance, the results obtained by SLBLR and the Column Generation method \cite{Sadykov15} are comparable. For instance d401600, SLBLR obtains a better feasible solution and for instance d801600, the advantage over the existing methods is even more pronounced. 

\begin{table}[!htbp]

\caption{Comparison Against the Best Results Currently Available.}
\centering 
\begin{tabular}{>{\centering\arraybackslash}p{11mm} | >{\centering\arraybackslash}p{17mm}|>{\centering\arraybackslash}p{14mm}|>{\centering\arraybackslash} p{17mm}|>{\centering\arraybackslash}p{17mm}|>{\centering\arraybackslash}p{17mm}}
\hline\hline 
 & New Method  & Posta \cite{Posta12} & Sadykov \cite{Sadykov15} & Haddadi \cite{Haddadi19}  &   Bragin \cite{Bragin19} \\
  [0.5ex] 

 Instance & (SLBLR)  & (DFLBnB) & (Column  & (VLSN) & (SAVLR) \\
&   &  & Generation) &  &  \\
\hline 
d201600  &     \textbf{97825} (1195)	& $-$$^{\dagger}$ & \textbf{97825} (1026) & 97836 (5364) & 97828 (1371)
    \\ 
d401600  &    \textbf{97105$^*$} (836) &  $-$$^{\dagger}$ &	97106 (919) & 97125 (5364) & 97111 (1183)
    \\ 
d801600 &  \textbf{97034$^*$} (3670) &	 $-$$^{\dagger}$ & 97037 (10860) & 97075 (5364) & 97039 (1350)
     \\ 
e201600  &     \textbf{180645$^{**}$} (85) & \textbf{180645} (40) & $-$ & \textbf{180645} (749)	& $-$ 
    \\ 
e401600  &    \textbf{178293$^{**}$} (2478) & \textbf{178293} (243) &	$-$ &	\textbf{178293} (749) & $-$ 
    \\ 
e801600 &     \textbf{176820$^{**}$} (1762) & \textbf{176820} (75) & $-$ & 176821 (749) &	$-$

     \\ 

\hline 
\end{tabular}
\label{table4} 
\begin{tablenotes}
   \item[$^*$]  The optimality is certified by the LP optimal values, which are 97105 and 97034 for instances d401600 and d801600, respectively.  
   \item[$^{**}$]  The optimality is certified through the lower bound results of, i.e., Posta et al. \cite[p. 160]{Posta12}. 
   \item[$-^{\dagger}$] Not solved to optimality within 24 hours and not reported within the original paper of \cite{Posta12}. 
   \item[$-$] These instances were not considered within \cite{Sadykov15} and \cite{Bragin19}. 
  \end{tablenotes}
\end{table}

To the best of the authors' knowledge, no other reported method obtained optimal results for instances d401600 and d801600. SLBLR outperforms SAVLR \cite{Bragin19} as well, thereby demonstrating that the fast convergence offered by the novel ``level-based'' stepsizing, with other things being equal, translates into better results as compared to those obtained by SAVLR, which employs the ``contraction mapping'' stepsizing \cite{Bragin19}. Lastly, the methods developed in \cite{Posta12, Haddadi19, Sadykov15} specifically target GAPs, whereas the SLBLR method developed in this paper has broader applicability.










\subsubsection{Stochastic Job-Shop Scheduling with the Consideration of Scrap and Rework} \label{JSS}
To demonstrate the computational capability of the method to solve large-scale stochastic MILP problems, a job-shop scheduling problem is considered. 
Within a job shop, each job requires a specific sequence of operations and the processing time for each operation. Operations are performed by a set of eligible machines.
To avoid late shipments,
expected tardiness is minimized.
Limited machine capacity brings a layer of difficulty since multiple ``individual-job'' subproblems are considered together competing for limited resources (machines). Another difficulty arises because of uncertainties, including processing times  \citep{Golenko02, lei11, zhang13, shen16, jamili19, horng21} and scrap \cite{Wilson22, Sun22, Bragin22}.
Re-manufacturing of one part may affect and disrupt the overall schedule within the entire job shop, thereby leading to unexpectedly high delays in production. 

\begin{figure}[!htbp]
		\centering
 \includegraphics[trim = 0 0 0 0, scale=.62]{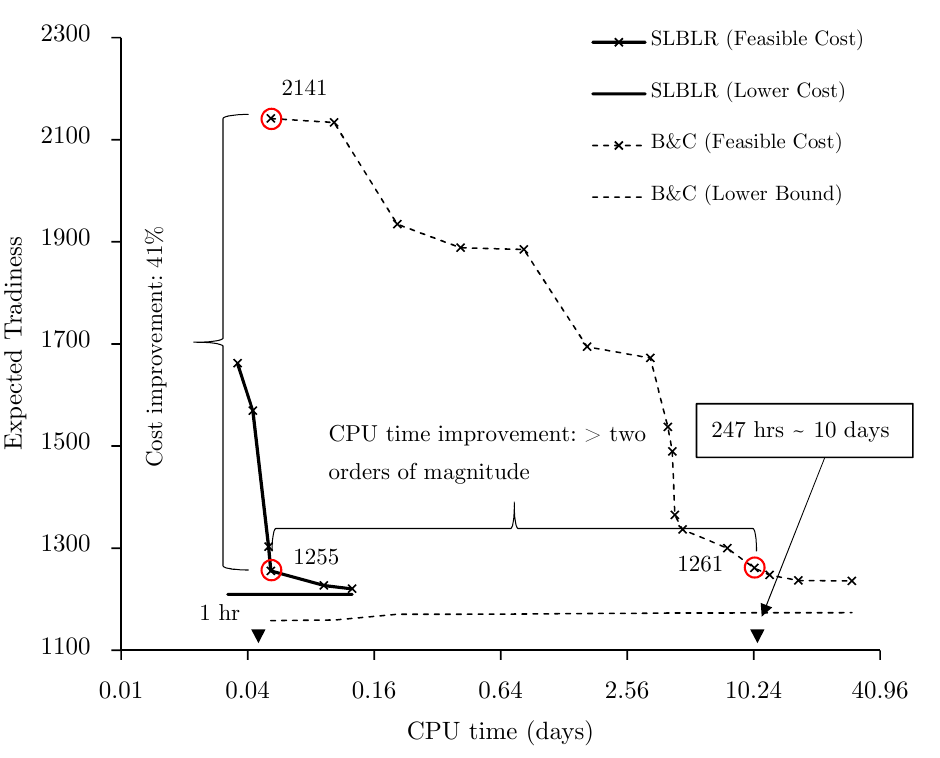}		
  \caption{The results for Example \ref{JSS} 
  are illustrated. SLBLR performs more than two orders of magnitude faster than branch-and-cut implemented in CPLEX.} 
  \label{fig:ex2_results}
\end{figure}

In this paper, we modified data from \cite{Luh93} by modifying several jobs by increasing the number of operations (e.g., from 1 to 6) and decreasing the capacities of a few machines; the data are in Tables \ref{duedatesex2} and \ref{proctimeex2}. The stochastic version of the problem with the consideration of scrap and rework is available at \cite{Bragin22}. With these changes,
the running time of CPLEX spans multiple days as demonstrated in Figure \ref{fig:ex2_results} 
In contrast, within the new method,
a solution of the same quality as that obtained by CPLEX, is obtained within roughly 1 hour of CPU time. The new method is operationalized by relaxing machine capacity constraints \cite[p. 14, eq. (16)]{Bragin22} and coordinating resulting job subproblems; at convergence, the beginning times of several jobs are adjusted by a few time periods to remove remaining machine capacity constraint violations. 

\subsubsection{Multi-Stage Pharmaceutical Scheduling}\label{pharma}

\begin{figure}[!htbp]
     \centering
     \begin{subfigure}[b]{0.0\textwidth}
         \includegraphics[trim = 0 0 0 0, scale=.4]{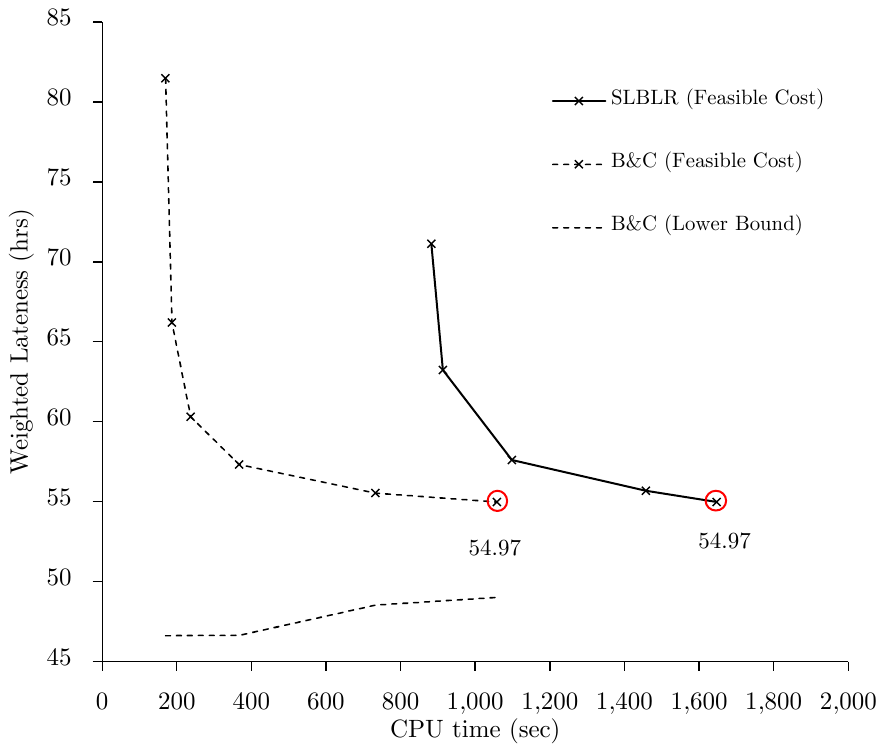}	      
     \end{subfigure}
     \hfill
     \begin{subfigure}[b]{0.5\textwidth}
         \includegraphics[trim = 0 0 0 0, scale=.4]{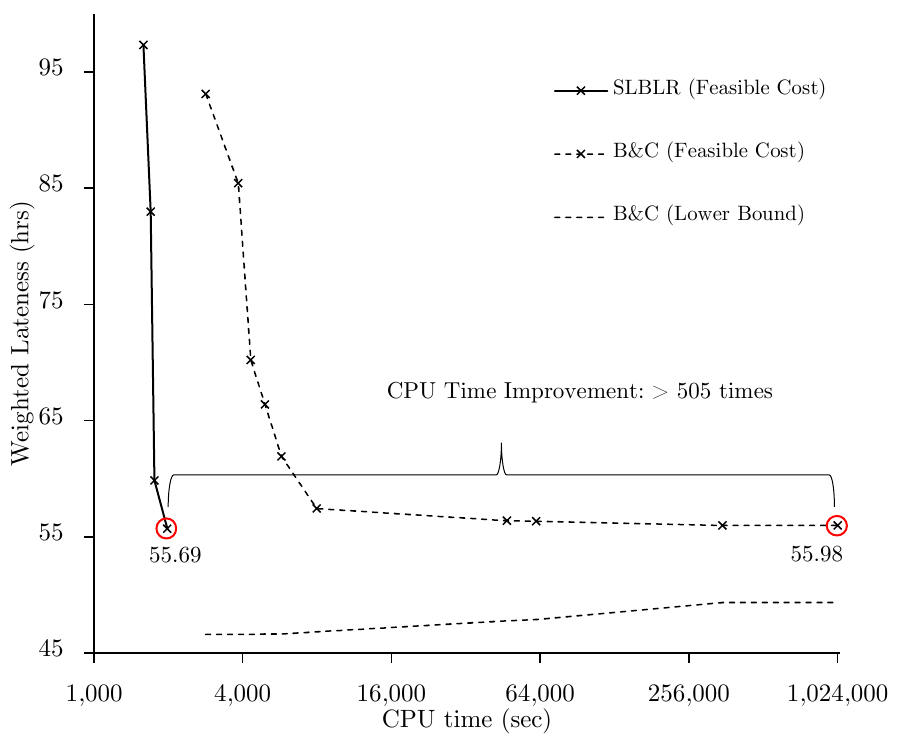}	
     \end{subfigure}
        \caption{The results for Example \ref{pharma} with 30 product orders (left) and 60 product orders (right) are illustrated. 
        }
        \label{fig:Fig4}
\end{figure}

To demonstrate the capability of the method to solve scheduling problems complicated by the presence of sequence-dependent setup times, a multi-stage pharmaceutical scheduling problem proposed by \cite{Kopanos10} is considered. 
Setup times vary based on the sequencing of products on each unit (machine). Scheduling in this context is combinatorial in the number of product orders, units, and stages. The new method is operationalized by relaxing constraints \cite[Eqs. (1)-(3), p. 646]{Kopanos10} that couple individual processing units, namely assignment, and processing time constraints. 
The formulation is available in \ref{Ex3_Formulation}, and the results with SLBLR and Branch-and-Cut are demonstrated in Figure \ref{fig:Fig4}. 

With a relatively small number of product orders, 30, an optimal solution with a feasible cost of 54.97 was found by CPLEX within 1057.78 seconds. The optimality is verified by running CPLEX until the gap is 0\%; it took 171993.27 seconds to verify the optimality. SLBLR takes a slightly longer time to obtain the same solution - 1647.35 seconds (Figure \ref{fig:Fig4} (left)). 
In contrast, with 60 product orders, CPLEX  no longer obtains good solutions in a computationally efficient manner; a solution with a feasible cost of 55.98 is obtained after 1,000,000 seconds. Within SLBLR, a solution with a feasible cost of 55.69 is obtained within 1978.04 seconds. This is more than two orders of magnitude of improvement over CPLEX. The drastic differences when the number of product orders doubles is demonstrated in Figure \ref{fig:Fig4} (right; log scale). When doubling the number of products, CPLEX's performance drastically deteriorated, while the performance of SLBLR is 
scalable.

\section{Discussion}\label{sec13}

This paper develops a novel MILP solution methodology based on the Lagrangian Relaxation method. Salient features of the novel SLBLR method, inherited from the previous versions of Lagrangian Relaxation, are 1. reduction of the computational effort required to obtain Lagrangian-multiplier-updating directions and 2. alleviation of zigzagging of multipliers. The key novel feature of the method, which the authors believe gives SLBLR the decisive advantage, is the innovative exploitation of the underlying geometric-convergence potential inherent to Polyak's stepsizing formula without the heuristic adjustment of hyperparameters for the estimate of $q(\lambda^*)$ - the associated the ``level'' values are determined purely through the simple auxiliary ``multiplier-divergence-detection'' constraint satisfaction problem. 
Through testing, it is discovered that SLBLR is robust with respect to the choice of initial stepsizes and multipliers, computationally efficient, competitive, and general. Several problems from diverse disciplines are tested and the superiority of SLBLR is demonstrated. 
While ``separable'' MILP problems are considered, no particular problem characteristics such as linearity or separability have been used to obtain ``level'' values, and thus SLBLR has the potential to solve a broad class of MIP problems.

\noindent \textbf{Acknowledgements} 
The work by M.A.B. was supported in part by the US NSF under award ECCS-1810108.

\noindent \textbf{Author Contributions} 
M.A.B. conceptualized the project, developed the novel methodology and conducted experiments. E.L.T. developed the fourth case study and supported testing and interpretation of results. M.A.B. drafted the initial manuscript with revisions from E.L.T. Both authors provided feedback on the final version. All authors read and approved the final manuscript.

\noindent \textbf{Competing Interests} 
The authors declare no competing interests.

\noindent \textbf{Data Availability.} Data supporting the results of Example \ref{Example:Ex_GAP} are located at \href{https://www-or.amp.i.kyoto-u.ac.jp/members/yagiura/gap/}{https://www-or.amp.i.kyoto-u.ac.jp/members/yagiura/gap/}; for Example \ref{JSS}, data are located in Tables \ref{duedatesex2} and \ref{proctimeex2} as well as in subsection \ref{53}; for Example \ref{pharma}, data are taken from \cite{Kopanos10}. CPLEX codes for Example \ref{SmallEx} as well as for the GAP instance d201600 
 from Example \ref{Example:Ex_GAP} are located at \href{https://github.com/ProfBragin/SLBLR}{https://github.com/ProfBragin/SLBLR } 

\bibliographystyle{naturemag}

\pagebreak









\newcommand{\beginsupplement}{%
        \setcounter{table}{0}
        \renewcommand{\thetable}{S\arabic{table}}%
        \setcounter{figure}{0}
        \renewcommand{\thefigure}{S\arabic{figure}}%
     }

\section{Supplementary Information.} \label{Supp}
\beginsupplement

\subsection{Methods: Previous Methods for Non-smooth Optimization and for MILP} \label{OnlineMethods}

\noindent \textbf{The 1990s: The Subgradient-Level Method.} The subgradient-level method, developed by \cite{Goffin99}, overcomes the unavailability of the knowledge about the optimal value needed to compute Polyak's stepsize \eqref{eq7} by adaptively adjusting a ``level'' estimate based on the detection of ``sufficient descent'' of the function and ``oscillation'' of solutions.  

In terms of the problem \eqref{eq3}, the procedure of the method is explained as follows: the ``level'' estimate $q^{k}_{lev} = q^{k_j}_{rec} + \delta_j$ is used in place of the optimal dual value $q(\lambda^{*})$, where $q^{k}_{rec}$ is the best dual value (``record objective value'') obtained up to an iteration $k,$ and $\delta_j$ is an adjustable parameter with $j$ denoting the $j^{th}$ update of $q^{k}_{lev}.$ The main premise behind this is when $\delta_j$ is ``too large,'' then multipliers will exhibit oscillations while traveling a significant (predefined) distance $R$ without improving the ``record'' value. In this case, the parameter $\delta_j$ is updated as $\delta_{j+1} = \beta \cdot \delta_j$ with $\beta = \frac{1}{2}.$ On the other hand, if $\delta_j$ is such that the dual value is sufficiently increased: $q(\lambda^k) \geq q^{k}_{lev} + \tau \cdot \delta_j,$ with $\tau = \frac{1}{2},$ then the parameter $\delta_j$ is unchanged and the distance traveled by multipliers is reset to 0 to avoid premature reduction of $\delta_j$ by $\beta$ in future iterations.

\noindent \textbf{The Early 2000s: Incremental Subgradient Methods \cite{Nedic01b, Nedic01}.}
The main idea of the incremental subgradient method is to improve convergence by solving a subproblem $i$ before updating multipliers.  After one subgradient component is updated, rather than updating all the multipliers ``at once,'' within the incremental subgradient methods, multipliers are updated ``incrementally.'' After the $i^{th}$ subgradient component is calculated, the multipliers are incrementally updated as 
\vspace{-2mm}
\begin{flalign}
& \psi^{k}_i = \psi^k_{i-1} + s^k \cdot \left(A_i x_i^k - \beta_i \right).   \label{eq11}
\end{flalign}
Here $\beta_i$ are the vectors such that $\sum_{i=1}^I \beta_i = b$, for example, $\beta_i = \frac{b}{I}.$  Only after all $i$ subproblems are solved, are the multipliers ``fully'' updated as  
\begin{flalign}
& \lambda^{k+1} = \psi^{k}_I.   \label{eq12a}
\end{flalign}
Convergence results of the subgradient-level method (\cite{Goffin99}) have been extended for the subgradient methods and proved. Variations of the method were proposed with $\beta$ and $\tau$ belonging to an interval $[0,1]$ rather than being equal to $\frac{1}{2}.$ Moreover, to improve convergence, rather than using constant $R,$ a sequence of $R_l$ was proposed such that $\sum_{l=1}^{\infty} R_l = \infty.$ While the method reduces the effort by solving one subproblem as a time, in order to compute the ``level'' values, the ``record'' dual value is required; in order to obtain dual values, all subproblems need to be solved optimally without updating the multipliers. 

\noindent \textbf{The 2010s: The Surrogate Lagrangian Relaxation Method \cite{SLR}.}
Convergence of the method is based on the ``contraction mapping'' concept. Namely, within the method, distances between multipliers at consecutive iterations are required to decrease, i.e., 
\begin{flalign}
& \| \lambda^{k+1} - \lambda^k \| = \alpha_k \cdot \| \lambda^{k} - \lambda^{k-1} \|, \quad 0 \leq \alpha_k < 1. \label{eq12b}
\end{flalign}
Based on \eqref{eq10}, the stepsizing formula has been derived: 
\begin{flalign}
& s^k = \alpha_k \cdot \frac{s^{k-1}\|g(\tilde{x}^{k-1})\|}{\|g(\tilde{x}^{k})\|}. \label{eq13a}
\end{flalign}
Moreover, a specific formula to set $\alpha_k$ has been developed to guarantee convergence: 
\begin{flalign} 
& \alpha_k = 1-\frac{1}{M \cdot k^{1-\frac{1}{k^r}}}, \; M \geq 1, \; 0 
\leq r \leq 1. \label{eq14}
\end{flalign}
Within \cite[Figs. 3-5,7; pp. 195-199]{SLR}, the SLR method showed the advantage against the above-mentioned subgradient-level method \cite{Goffin99} and the incremental subgradient methods \cite{Nedic01b, Nedic01} for non-smooth optimization. Unlike the methods of \cite{Nedic01b, Nedic01}, the SLR method does not require obtaining dual values to set stepsizes, which further reduces the effort. 

\noindent \textbf{Surrogate Absolute-Value Lagrangian Relaxation \cite{Bragin19}.} Aiming to simultaneously guarantee convergence while ensuring a fast reduction of constraint violations and preserving the linearity of the original MILP problem, the Surrogate Absolute-Value Lagrangian Relaxation (SAVLR) method was developed. Within the method, the following dual problem is considered:
\begin{flalign}
& \max_{\lambda} \{q_{\rho}(\lambda): \lambda \in \Omega \subset \mathbb{R}^m\},  \label{eq122}
\end{flalign}
where
\begin{flalign}
& q_{\rho}(\lambda) = 
\min_{(x,y)} \Bigg\{\sum_{i=1}^I \left((c_i^x)^T x_i + (c_i^y)^T y_i\right) +
\nonumber  \lambda^T \cdot \left(\sum_{i=1}^I A_i^x x_i + \sum_{i=1}^I A_i^y y_i - b\right) +
 \\ & \rho \cdot \| \sum_{i=1}^I A_i^x x_i + \sum_{i=1}^I A_i^y y_i - b \|_1, \left\{x_i,y_i\right\} \in \mathcal{F}_i, i = 1,\dots,I \Bigg\}.  \label{eq123}
\end{flalign}

The above minimization involves $l_1$-absolute-value piece-wise linear penalties, which efficiently penalize constraint violations and are exactly linearizable thereby enabling the use of MILP solvers. Within \cite[Figs. 1, 2, 5-8; pp. 537-539]{Bragin19}, the SAVLR method showed the advantage over the above SLR method to solve GAPs. For completeness of comparison with the newly-developed SLBLR method, SAVLR is also included in Table \ref{table4}. 

\subsection{Example \ref{Example:Ex_GAP} Description and Tables.}\label{Formulation_GAP}

The mathematical formulation of GAP is: 
\begin{flalign} 
& \min_{x_{i,j}} \sum_{i=1}^I \sum_{j=1}^J g_{i,j} x_{i,j}, \;\; x_{i,j} \in \{0,1\}, \;\; g_{i,j} \geq 0, \label{eq35a}
\end{flalign}
\begin{flalign} 
& s.t. \;\; \sum_{i=1}^I a_{i,j} x_{i,j} \leq b_j, \;\; j = 1,\dots, J, \;\; a_{i,j} \geq 0, \; b_j \geq 0, \label{eq35b}
\end{flalign}
\begin{flalign} 
& \sum_{j=1}^J x_{i,j} = 1, \;\; j = 1,\dots, J. \label{eq35c}
\end{flalign}
The objective \eqref{eq35a} is to minimize the total assignment cost by deciding which job $i$ is to be assigned to which machine $j$; if job $i$ is assigned to machine $j$, then $x_{i,j} = 1$ and $x_{i,j} = 0$, otherwise. Constraints \eqref{eq35b} ensure that the total amount of time required by all jobs to be processed on machine $j$ should not exceed the total machine's available time $b_j.$ The assignment constraints \eqref{eq35c} ensure that each job is assigned to only one machine. Within SLBLR, constraints \eqref{eq35c} are relaxed. 

\subsection{Example \ref{JSS} Description and Tables.} \label{53}
For the job-shop scheduling example, to create a difficult testing case, data from \cite{Luh93} are considered and several parameters are modified. Namely, 16 jobs: Jobs 1 through 6 and 118 through 127 originally with one operation, are modified by adding 5 more operations, with processing times generated randomly from a range $[1,5]$ based on a discrete uniform distribution. Moreover, to create a case mimicking ``labor shortage'', the machine capacity of certain machine types is reduced by 2, and the resulting machine capacities are $M = \{
1,	1,	\textbf{2},	2,	1,	1,	\textbf{2},	1,	1,	\textbf{4},	1,	1,	1,	1,	2,	1,	1\},$ with modified capacities shown in bold. The due dates and the number of operations per job are shown in Table \ref{duedatesex2}. Several due dates, originally negative, are modified. Other parameters such as processing times are shown in Table \ref{proctimeex2}.

\begin{table}[!htbp]
\caption{Data for Example \ref{JSS}.}
\label{duedatesex2}

\centering
		{\def\arraystretch{0.4}  
\begin{tabular*}{1.05\textwidth}{@{\extracolsep{\fill}}ccc| ccc| ccc| ccc}
\hline
\hline
$i$ & $d_i$ & $J_i$ & $i$ & $d_i$ & $J_i$ & $i$ & $d_i$ & $J_i$ & $i$ & $d_i$ & $J_i$ \\ 
\hline
1	&	10	&	6	&	33	&	37	&	2	&	65	&	-81	&	1	&	97	&	14	&	1	\\
2	&	10	&	6	&	34	&	37	&	2	&	66	&	37	&	1	&	98	&	14	&	1	\\
3	&	9	&	6	&	35	&	159	&	2	&	67	&	37	&	1	&	99	&	37	&	1	\\
4	&	61	&	6	&	36	&	8	&	1	&	68	&	71	&	2	&	100	&	-1	&	1	\\
5	&	61	&	6	&	37	&	12	&	1	&	69	&	88	&	2	&	101	&	2	&	1	\\
6	&	-3	&	6	&	38	&	20	&	1	&	70	&	18	&	1	&	102	&	-4	&	1	\\
7	&	-85	&	3	&	39	&	40	&	1	&	71	&	8	&	2	&	103	&	4	&	1	\\
8	&	10	&	1	&	40	&	62	&	1	&	72	&	31	&	2	&	104	&	7	&	1	\\
9	&	18	&	1	&	41	&	84	&	1	&	73	&	35	&	2	&	105	&	53	&	2	\\
10	&	145	&	3	&	42	&	35	&	1	&	74	&	36	&	1	&	106	&	95	&	2	\\
11	&	145	&	3	&	43	&	35	&	1	&	75	&	39	&	2	&	107	&	132	&	2	\\
12	&	-169	&	2	&	44	&	157	&	1	&	76	&	40	&	2	&	108	&	174	&	2	\\
13	&	-114	&	1	&	45	&	157	&	1	&	77	&	42	&	2	&	109	&	218	&	2	\\
14	&	-123	&	1	&	46	&	157	&	1	&	78	&	44	&	2	&	110	&	261	&	2	\\
15	&	91	&	3	&	47	&	157	&	1	&	79	&	31	&	2	&	111	&	304	&	2	\\
16	&	91	&	3	&	48	&	35	&	1	&	80	&	35	&	2	&	112	&	28	&	2	\\
17	&	91	&	3	&	49	&	35	&	1	&	81	&	36	&	2	&	113	&	2	&	1	\\
18	&	91	&	3	&	50	&	7	&	1	&	82	&	39	&	2	&	114	&	2	&	1	\\
19	&	91	&	3	&	51	&	-53	&	1	&	83	&	40	&	2	&	115	&	111	&	1	\\
20	&	91	&	3	&	52	&	70	&	1	&	84	&	42	&	2	&	116	&	111	&	1	\\
21	&	-65	&	3	&	53	&	3	&	1	&	85	&	44	&	2	&	117	&	-6	&	2	\\
22	&	-51	&	3	&	54	&	-4	&	1	&	86	&	71	&	1	&	118	&	17	&	6	\\
23	&	-47	&	1	&	55	&	12	&	1	&	87	&	71	&	1	&	119	&	-1	&	6	\\
24	&	-35	&	1	&	56	&	29	&	2	&	88	&	88	&	1	&	120	&	-27	&	6	\\
25	&	2	&	1	&	57	&	132	&	1	&	89	&	88	&	1	&	121	&	-6	&	6	\\
26	&	6	&	1	&	58	&	-22	&	1	&	90	&	71	&	2	&	122	&	-6	&	6	\\
27	&	-47	&	2	&	59	&	13	&	1	&	91	&	81	&	2	&	123	&	1	&	6	\\
28	&	38	&	1	&	60	&	-10	&	1	&	92	&	11	&	1	&	124	&	56	&	6	\\
29	&	38	&	1	&	61	&	12	&	1	&	93	&	8	&	1	&	125	&	78	&	6	\\
30	&	161	&	1	&	62	&	12	&	1	&	94	&	8	&	1	&	126	&	98	&	6	\\
31	&	-45	&	1	&	63	&	6	&	6	&	95	&	11	&	1	&	127	&	19	&	6	\\
32	&	-53	&	1	&	64	&	35	&	1	&	96	&	8	&	1	&		&		&		\\

\hline
\hline
\end{tabular*}
\label{table:mechanical_property}
}
\end{table}

\begin{table}[!htbp]
\caption{Data for Example \ref{JSS}.}
\label{proctimeex2}
\centering
		{\def\arraystretch{0.75}  
\begin{tabular*}{1\textwidth}{@{\extracolsep{\fill}}cccc| cccc| cccc| cccc}
\hline
\hline
$i$ & $j$ & $p_{i,j}$ & $m$ & $i$ & $j$ & $p_{i,j}$ & $m$ & $i$ & $j$ & $p_{i,j}$ & $m$ & $i$ & $j$ & $p_{i,j}$ & $m$ \\
\hline
1	&	1	&	1	&	1	&	20	&	2	&	1	&	2	&	69	&	2	&	1	&	10	&	112	&	1	&	2	&	5	\\
1	&	2	&	3	&	1	&	20	&	3	&	1	&	3	&	70	&	1	&	1	&	8	&	112	&	2	&	1	&	10	\\
1	&	3	&	3	&	1	&	21	&	1	&	1	&	2	&	71	&	1	&	3	&	15	&	113	&	1	&	5	&	10	\\
1	&	4	&	3	&	1	&	21	&	2	&	1	&	3	&	71	&	2	&	1	&	1	&	114	&	1	&	3	&	10	\\
1	&	5	&	4	&	1	&	21	&	3	&	2	&	2	&	72	&	1	&	1	&	13	&	115	&	1	&	4	&	2	\\
1	&	6	&	1	&	1	&	22	&	1	&	1	&	3	&	72	&	2	&	1	&	13	&	116	&	1	&	3	&	3	\\
2	&	1	&	6	&	1	&	22	&	2	&	1	&	2	&	73	&	1	&	1	&	13	&	117	&	1	&	1	&	9	\\
2	&	2	&	5	&	1	&	22	&	3	&	2	&	3	&	73	&	2	&	1	&	13	&	117	&	2	&	1	&	9	\\
2	&	3	&	1	&	1	&	23	&	1	&	1	&	2	&	74	&	1	&	1	&	13	&	118	&	1	&	1	&	10	\\
2	&	4	&	2	&	1	&	24	&	1	&	1	&	11	&	75	&	1	&	1	&	13	&	118	&	2	&	6	&	10	\\
2	&	5	&	2	&	1	&	25	&	1	&	2	&	12	&	75	&	2	&	1	&	13	&	118	&	3	&	5	&	10	\\
2	&	6	&	6	&	1	&	26	&	1	&	2	&	2	&	76	&	1	&	1	&	13	&	118	&	4	&	4	&	10	\\
3	&	1	&	2	&	1	&	27	&	1	&	28	&	8	&	76	&	2	&	1	&	13	&	118	&	5	&	1	&	10	\\
3	&	2	&	1	&	1	&	27	&	2	&	3	&	9	&	77	&	1	&	1	&	13	&	118	&	6	&	2	&	10	\\
3	&	3	&	1	&	1	&	28	&	1	&	1	&	1	&	77	&	2	&	1	&	13	&	119	&	1	&	2	&	10	\\
3	&	4	&	3	&	1	&	29	&	1	&	1	&	1	&	78	&	1	&	1	&	13	&	119	&	2	&	4	&	10	\\
3	&	5	&	5	&	1	&	30	&	1	&	1	&	1	&	78	&	2	&	1	&	13	&	119	&	3	&	2	&	10	\\
3	&	6	&	5	&	1	&	31	&	1	&	1	&	8	&	79	&	1	&	1	&	13	&	119	&	4	&	5	&	10	\\
4	&	1	&	1	&	2	&	32	&	1	&	10	&	5	&	79	&	2	&	1	&	13	&	119	&	5	&	3	&	10	\\
4	&	2	&	5	&	2	&	33	&	1	&	2	&	10	&	80	&	1	&	1	&	13	&	119	&	6	&	2	&	10	\\
4	&	3	&	3	&	2	&	33	&	2	&	2	&	14	&	80	&	2	&	1	&	13	&	120	&	1	&	2	&	2	\\
4	&	4	&	2	&	2	&	34	&	1	&	2	&	10	&	81	&	1	&	1	&	13	&	120	&	2	&	3	&	2	\\
4	&	5	&	1	&	2	&	34	&	2	&	2	&	14	&	81	&	2	&	1	&	13	&	120	&	3	&	5	&	2	\\
4	&	6	&	4	&	2	&	35	&	1	&	2	&	10	&	82	&	1	&	1	&	13	&	120	&	4	&	1	&	2	\\
5	&	1	&	6	&	3	&	35	&	2	&	2	&	14	&	82	&	2	&	1	&	13	&	120	&	5	&	2	&	2	\\
5	&	2	&	6	&	3	&	36	&	1	&	14	&	10	&	83	&	1	&	1	&	13	&	120	&	6	&	3	&	2	\\
5	&	3	&	3	&	3	&	37	&	1	&	7	&	10	&	83	&	2	&	1	&	13	&	121	&	1	&	3	&	3	\\
5	&	4	&	5	&	3	&	38	&	1	&	7	&	10	&	84	&	1	&	1	&	13	&	121	&	2	&	3	&	3	\\
5	&	5	&	2	&	3	&	39	&	1	&	7	&	10	&	84	&	2	&	1	&	13	&	121	&	3	&	1	&	3	\\
5	&	6	&	3	&	3	&	40	&	1	&	7	&	10	&	85	&	1	&	1	&	13	&	121	&	4	&	6	&	3	\\
6	&	1	&	3	&	2	&	41	&	1	&	7	&	10	&	85	&	2	&	1	&	13	&	121	&	5	&	4	&	3	\\
6	&	2	&	1	&	2	&	42	&	1	&	2	&	10	&	86	&	1	&	1	&	10	&	121	&	6	&	4	&	3	\\
6	&	3	&	6	&	2	&	43	&	1	&	2	&	10	&	87	&	1	&	1	&	10	&	122	&	1	&	1	&	2	\\
6	&	4	&	5	&	2	&	44	&	1	&	2	&	10	&	88	&	1	&	1	&	10	&	122	&	2	&	1	&	2	\\
6	&	5	&	1	&	2	&	45	&	1	&	2	&	10	&	89	&	1	&	1	&	10	&	122	&	3	&	2	&	2	\\
6	&	6	&	6	&	2	&	46	&	1	&	2	&	10	&	90	&	1	&	5	&	7	&	122	&	4	&	5	&	2	\\
7	&	1	&	2	&	5	&	47	&	1	&	2	&	10	&	90	&	2	&	1	&	10	&	122	&	5	&	2	&	2	\\
7	&	2	&	4	&	6	&	48	&	1	&	2	&	10	&	91	&	1	&	5	&	7	&	122	&	6	&	2	&	2	\\
7	&	3	&	2	&	6	&	49	&	1	&	2	&	10	&	91	&	2	&	1	&	10	&	123	&	1	&	1	&	3	\\
8	&	1	&	1	&	7	&	50	&	1	&	2	&	10	&	92	&	1	&	2	&	10	&	123	&	2	&	4	&	3	\\
9	&	1	&	6	&	8	&	51	&	1	&	1	&	2	&	93	&	1	&	2	&	10	&	123	&	3	&	3	&	3	\\
10	&	1	&	4	&	8	&	52	&	1	&	7	&	3	&	94	&	1	&	2	&	10	&	123	&	4	&	2	&	3	\\
10	&	2	&	2	&	6	&	53	&	1	&	2	&	10	&	95	&	1	&	2	&	10	&	123	&	5	&	4	&	3	\\
10	&	3	&	1	&	8	&	54	&	1	&	2	&	10	&	96	&	1	&	2	&	10	&	123	&	6	&	6	&	3	\\
11	&	1	&	4	&	8	&	55	&	1	&	1	&	10	&	97	&	1	&	1	&	10	&	124	&	1	&	1	&	9	\\
11	&	2	&	2	&	6	&	56	&	1	&	4	&	5	&	98	&	1	&	1	&	10	&	124	&	2	&	4	&	9	\\
11	&	3	&	1	&	8	&	56	&	2	&	3	&	10	&	99	&	1	&	4	&	2	&	124	&	3	&	3	&	9	\\
12	&	1	&	28	&	8	&	57	&	1	&	13	&	16	&	100	&	1	&	2	&	3	&	124	&	4	&	2	&	9	\\
12	&	2	&	3	&	9	&	58	&	1	&	10	&	5	&	101	&	1	&	2	&	2	&	124	&	5	&	4	&	9	\\
13	&	1	&	1	&	10	&	59	&	1	&	13	&	7	&	102	&	1	&	1	&	3	&	124	&	6	&	6	&	9	\\
14	&	1	&	1	&	10	&	60	&	1	&	3	&	4	&	103	&	1	&	2	&	2	&	125	&	1	&	1	&	9	\\
15	&	1	&	1	&	3	&	61	&	1	&	13	&	4	&	104	&	1	&	2	&	3	&	125	&	2	&	4	&	9	\\
15	&	3	&	1	&	2	&	62	&	1	&	13	&	4	&	105	&	1	&	1	&	10	&	125	&	3	&	3	&	9	\\
16	&	1	&	1	&	3	&	63	&	1	&	1	&	1	&	105	&	2	&	4	&	5	&	125	&	4	&	2	&	9	\\
16	&	2	&	1	&	2	&	63	&	2	&	1	&	1	&	106	&	1	&	1	&	10	&	125	&	5	&	4	&	9	\\
16	&	3	&	1	&	3	&	63	&	3	&	1	&	1	&	106	&	2	&	4	&	5	&	125	&	6	&	6	&	9	\\
17	&	1	&	1	&	2	&	63	&	4	&	1	&	1	&	107	&	1	&	1	&	10	&	126	&	1	&	1	&	9	\\
17	&	2	&	1	&	3	&	63	&	5	&	1	&	1	&	107	&	2	&	4	&	5	&	126	&	2	&	4	&	9	\\
17	&	3	&	1	&	2	&	63	&	6	&	1	&	1	&	108	&	1	&	1	&	10	&	126	&	3	&	3	&	9	\\
18	&	1	&	1	&	3	&	64	&	1	&	2	&	9	&	108	&	2	&	4	&	5	&	126	&	4	&	2	&	9	\\
18	&	2	&	1	&	2	&	65	&	1	&	1	&	10	&	109	&	1	&	1	&	10	&	126	&	5	&	4	&	9	\\
18	&	3	&	1	&	3	&	66	&	1	&	1	&	10	&	109	&	2	&	4	&	5	&	126	&	6	&	6	&	9	\\
19	&	1	&	1	&	2	&	67	&	1	&	1	&	10	&	110	&	1	&	1	&	10	&	127	&	1	&	1	&	17	\\
19	&	2	&	1	&	3	&	68	&	1	&	5	&	7	&	110	&	2	&	4	&	5	&	127	&	2	&	4	&	17	\\
19	&	3	&	1	&	2	&	68	&	2	&	1	&	10	&	111	&	1	&	1	&	10	&	127	&	3	&	3	&	17	\\
20	&	1	&	1	&	3	&	69	&	1	&	5	&	7	&	111	&	2	&	4	&	5	&	127	&	4	&	2	&	17	\\

\hline
\hline
\end{tabular*}
}
\end{table}

\subsection{Example \ref{pharma} Description.}\label{Ex3_Formulation}

The data and formulation for the case example of multi-stage pharmaceutical scheduling, presented in Section \ref{pharma}, come from \cite{Kopanos10}. In the following, we use terminology and notation consistent with its original presentation. A set of product orders $I$ are manufactured across a set of $S$ stages of production. At each stage, product $i\in I$ is processed on exactly one unit $j\in J$, representing a machine. Each order $i \in I$ has a due date, $\delta_i$. Deviations from this due date are penalized with a unit weight of $w^+$ (tardiness) and $w^-$ (earliness). Setup times for each order have sequence-independent and sequence-dependent components. Sequence-independent setup times, $\pi_i$, are only based on the order $i\in i$. The sequence-dependent setup times, $\gamma_{s,i,i'}$, for order $i'\in I$ in stage $s\in S$ vary based on the order $i\in I$ it follows. The processing time of order $i \in I$ on unit $j \in J$ is given by $\tau_{i,j}$. A Big-M parameter $M$ is defined as a sufficiently large number.

Binary assignment variables, $y_{j,i,s}$, indicate if unit $j\in J$ processes order $i\in I$ in stage $s\in S$. Binary precedence variables, $x_{j,i,i'}$, designate if order $i$ is processed before $i'\in I$ when both $i,i'$ are assigned to unit $j\in J$. The time that order $i\in I$ completes stage $s\in S$ is given by continuous variables, $c_{is}$. The tardiness and earliness of each order $i\in I$ are given by continuous variables, $t_i$ and $e_i$, respectively.

\begin{flalign} 
& \min_{t_i,e_i,y_{j,i,s},x_{j,i,i'}} \sum_{i \in I} \{w^+ t_i + w^- e_i\}
, \label{OF_pharma}
\end{flalign}

\begin{flalign} 
& s.t. \;\; e_i \geq \delta_i - c_{i,S_i^{last}}, \; \; \forall i\in I, \label{con_early}
\end{flalign}

\begin{flalign}
 t_i \geq c_{i,S_i^{last}} - \delta_i, \; \; \forall i\in I, \label{con_late}
\end{flalign}

\begin{flalign}
\sum_{j\in J} y_{j,i,s} = 1, \; \; \forall i \in I, s\in S, \label{con_assignorder}
\end{flalign}

\begin{flalign}
c_{i,s} \geq \sum_{j\in J} (\tau_{i,j} + \pi_i)y_{j,i,s}, \; \; \forall i\in i, s=1, \label{con_ctime1}
\end{flalign}

\begin{flalign}
c_{i,s} - \sum_{j\in J} (\tau_{i,j} + \pi_i)y_{j,i,s} = c_{i,s-1}, \; \; \forall i\in i, s\in S\setminus \{1\}, \label{con_ctimeafter1}
\end{flalign}

\begin{multline}
c_{i,s} + \gamma_{s,i,i'} \leq c_{i',s} - \pi_{i'} - \tau_{i',j} + M(1-x_{j,i,i'}) + M(2-y_{j,i,s} -y_{j,i',s}), \\
\; \; \forall s\in S, j\in J, i, i' \in I: i'\geq i+1, \label{con_prec1}
\end{multline}

\begin{multline}
c_{i',s} + \gamma_{s,i',i} \leq c_{i,s} - \pi_{i} - \tau_{i,j} + M(1-x_{j,i,i'}) + M(2-y_{j,i,s} -y_{j,i',s}), \\
\; \; \forall s\in S, j\in J, i, i' \in I: i'\geq i+1, \label{con_prec2}
\end{multline}

\begin{flalign}
y_{j,i,s},x_{j,i,i'} \in \{0,1\}, \; \; \forall j\in J, i,i'\in I, s\in S, \; \; \label{con_domain1}
\end{flalign}

\begin{flalign}
t_i, e_i, c_{i,s} \geq 0, \; \; \forall i\in I, s\in S. \label{con_domain2}
\end{flalign}

The objective function \eqref{OF_pharma} minimizes the total weighted deviations of orders from their due dates (referred to as in the original paper as ``weighted lateness''). Constraints \eqref{con_early} and \eqref{con_late} record the earliness and tardiness, respectively, of each order $i\in I$. Constraints \eqref{con_assignorder} require each order $i \in I$ to be assigned exactly one unit $j \in J$ in each stage $s\in S$. Constraints \eqref{con_ctime1} ensure the first stage completion times of each order $i\in I$ are at least the sum of the independent setup and processing times on the assigned unit $j\in J$. Similar constraints \eqref{con_ctimeafter1} ensure subsequent stage completion times as a function of independent setup, processing, and completion times of the previous stage. Constraints \eqref{con_prec1}-\eqref{con_prec2} enforce the sequencing between orders $i,i'\in I$ for each unit $j\in J$ and each stage $s\in S$. Constraints \eqref{con_domain1}-\eqref{con_domain2} are standard domain constraints.

There are 17 units and 6 stages of production. The problem is run with both 30 and 60 product orders. Minor changes from \cite{Kopanos10} are the incorporation of sequence-independent setup times and a non-optional stage three. We incorporated sequence-independent setup times by selecting the sequence-dependent setup times for one product in one stage, i.e., P01 in stage five. We generated sequence-dependent, stage three processing times for orders $i\in I$ that skipped stage three in \cite{Kopanos10} by randomly sampling from the two values given for the other products. To enforce the unallowable product-unit combinations from \cite{Kopanos10}, we set sufficiently large processing times where appropriate. 

\subsection{Proof of Corollary 1.} 

\proof{} 
Define two predicates
\begin{flalign} 
& A = \left\{s^k < \gamma \cdot \frac{q(\lambda^{*}) - L(\tilde{x}^k,\tilde{y}^k,\lambda^k)}{\|g(\tilde{x}^k,\tilde{y}^k)\|^2}\right\}.\label{eq21}
\end{flalign}
and 
\begin{flalign} 
& B = \left\{\|\lambda^{*}-\lambda^{k+1}\| < \|\lambda^{*}-\lambda^{k}\|\right\}.\label{eq22}
\end{flalign}
From Theorem 1 the following is true $A \Rightarrow B$ and from Corollary 1 the following is true $\neg B \Rightarrow \neg A.$ It remains to prove that both assertions are equivalent. Taking negation of $A \Rightarrow B$ leads to $\neg A \lor B,$ and taking negation of $\neg B \Rightarrow \neg A$ leads to $\neg (\neg B) \lor \neg A,$ which simplifies to $B \lor \neg A.$ 
\endproof

\end{document}